\documentclass[a4paper, 10pt,reqno]{amsart}  

\usepackage{amsmath} 
\usepackage{amssymb}
\usepackage{amsthm}  
\usepackage{color}
\usepackage{url}
\usepackage{enumitem}
\usepackage{nameref}
\usepackage{hyperref}
\usepackage{mathabx} 

\newtheorem{theorem}{Theorem}
\newtheorem{lemma}[theorem]{Lemma}

\newtheorem{proposition}[theorem]{Proposition}
\newtheorem{corollary}[theorem]{Corollary}
\newtheorem{definition}[theorem]{Definition}

\newtheorem*{mainresult}{Main Result}

\newtheorem{example}[theorem]{Example}
\newtheorem{remark}[theorem]{Remark}

\newcommand{\calK}{\mathcal{K}}

\newcommand{\Rpn}{\R_{0}^{+}}
\newcommand{\R}{\mathbb{R}}

\newcommand{\C}{\mathbb{C}}

\newcommand{\e}{\mathrm{e}}

\DeclareMathOperator{\ran}{ran}

\DeclareMathOperator*{\esssupp}{ess\,supp}

\allowdisplaybreaks

\numberwithin{theorem}{section}

\begin{document}

\title{
Integral input-to-state stability of unbounded bilinear control systems
}

\thanks{This work has been supported by the German Research Foundation (DFG) via the joint grant JA 735/18-1 / SCHW 2022/2-1.}

\author[R.~Hosfeld]{Ren\'e Hosfeld}
\address[RH]{University of Wuppertal, School of Mathematics and Natural Siences, IMACM, Gau\ss \-str.\ 20, 42119 Wuppertal, Germany and Department of Mathematics, University of Hamburg, Bundesstra\ss e 55, 20146 Hamburg, Germany}
\email{hosfeld@uni-wuppertal.de}
\author[B.~Jacob]{Birgit Jacob}
\address[BJ]{University of Wuppertal, School of Mathematics and Natural Sciences, IMACM, Gau\ss \-str.\ 20, 42119 Wuppertal, Germany}
\email{bjacob@uni-wuppertal.de}
\author[F.L.~Schwenninger]{Felix L.~Schwenninger}
\address[FLS]{ Department of Applied Mathematics, University of Twente, P.O. Box 217, 7500 AE Enschede, The Netherlands and Department of Mathematics, University of Hamburg, Bundesstra\ss e 55, 20146 Hamburg, Germany}
\email{f.l.schwenninger@utwente.nl}

\begin{abstract}
We study integral input-to-state stability of bilinear systems with  unbounded control operators and derive natural sufficient conditions. The results are applied to a bilinearly controlled Fokker-Planck equation.
\end{abstract} 




\begin{abstract}
We study integral input-to-state stability of bilinear systems with  unbounded control operators and derive natural sufficient conditions. The results are applied to a bilinearly controlled Fokker-Planck equation.
\end{abstract} 

\maketitle
\thispagestyle{plain}
\pagestyle{plain}




\section{Introduction}

%
 In this note we continue recent developments on {\em input-to-state stability (ISS)} for systems governed by evolution equations. This concept unifies both asymptotic stability with respect to the initial values and robustness with respect to the external inputs such as controls or disturbances. Loosely, if a system $\Sigma$ is viewed as a mapping which sends initial values $x_{0}\in X$ and inputs $u:[0,\infty)\to U$ to the time evolution $x:[0,T)\to X$ for some maximal $T>0$, then $\Sigma$ is ISS if $T=\infty$ and for all $t\in [0,\infty)$,
  \begin{equation}\label{intro:eq0}
  \|x(t)\|_{X} \leq \beta(\|x_{0}\|_{X},t)+\gamma(\sup_{s\in[0,t]}\|u(s)\|_{U}), \qquad \forall x_{0},u, 
  \end{equation}
  where  the continuous functions $\beta: \R_{0}^{+}\times\R_{0}^{+}\to \R_{0}^{+} $ and $\gamma:\R_{0}^{+}\to \R_{0}^{+}$ are of Lyapunov class $\mathcal{KL}$ and $\mathcal{K}$ respectively\footnote{I.e.\ 
   $\beta(0,s)=\lim_{t\to\infty}\beta(r,t)=\gamma(0)=0$ and $\beta(\cdot,t)$, $\beta\left(r,\frac{1}{\cdot}\right)$, $\gamma$ are strictly increasing on $\R^{+}$ for all $r,s,t>0$.
  }. Here $X$ is called the state space and $U$ the input space  equipped with norms $\|\cdot\|_{X}$ and $\|\cdot\|_{U}$. 
  For linear systems 
    \begin{equation*}\label{intro:eq1}
  \dot{x}(t)=Ax(t)+Bu(t),
  \end{equation*}
  where $A$ is the infinitesimal generator of a $C_{0}$-semigroup $(T(t))_{t\ge 0}$ on a Banach space $X$ and $B:U\to X$ is a bounded linear operator, 
  ISS is equivalent  to  uniform exponential stability of the semigroup \cite{DaM13,JNPS16}.  
If  $B$ is not  bounded as operator   from $U$ to $X$, which is  typically the case for boundary controlled PDEs, the property of being ISS becomes non-trivial even for linear systems. In  fact, this is closely related to suitable solution concepts see e.g.\ \cite{JNPS16,LhacShor19,Schw20}. Along with the recent developments in ISS theory for infinite-dimensional systems \cite{DaM13,DaM13b,GuivLogeOpme19,JaLoRy08,MiroWirt18a}, several partial results have been derived in the (semi)linear context, with a slight focus on parabolic equations, see e.g.\ \cite{JacoSchwZwar19,Krstic16MTNS,MaPr11,MiI15b,MiroKaraKrst17,ZhZh18}. We refer to recent  surveys on ISS for infinite-dimensional systems \cite{MiroPrie20,Schw20} and the book \cite{KaraKrst19book}.
The origin of ISS theory, introduced by Sontag in 1989 \cite{Sontag89ISS}, are non-linear systems and we refer the reader to \cite{sontagISS} for a survey on ISS for ODEs. 
Already seemingly harmless system classes such as bilinear systems 
\begin{equation}\label{intro:eq3}
  \dot{x}(t)=Ax(t)+B(x(t),u(t)),
\end{equation}
where $B(x,u)=\sum_{i=1}^{m}u_{i}B_{i}x$ and $A, B_{i}\in \R^{d\times d}$, see \cite{Elli09}, are typical counterexamples for ISS \cite{Son98}. Nevertheless, the following variant of ISS \cite{Son98} is satisfied by such systems; there exists functions $\beta\in \mathcal{KL}$ and $\gamma_{1},\gamma_{2}\in \mathcal{K}$ such that 
\begin{equation}\label{intro:eq5}
\|x(t)\|_{X} \leq \beta(\|x_{0}\|_{X},t)+\gamma_{1}\left(\int_{0}^{t}\gamma_{2}(\|u(s)\|_{U})\, \mathrm{d}s\right), \qquad \forall t>0,  \forall u,x_{0}, 
\end{equation}
which is called {\em integral input-to-state stable} (integral ISS), see also \cite{Son98}. Note that the terms involving $u$ in \eqref{intro:eq0} and \eqref{intro:eq5} cannot be compared for arbitrary $t>0$, general functions $u$,  and fixed functions $\gamma,\gamma_{1}, \gamma_{2}$. Still integral ISS and ISS are equivalent for infinite-dimensional linear systems with a bounded linear operator $B:U\to X$, \cite{JNPS16,MiI14b} as this reduces to uniform exponential stability of the uncontrolled system.  The corresponding question for general infinite-dimensional systems seems to be much harder and notorious questions remain, see \cite{JNPS16,NabSchwe,JacoSchwZwar19} and  \cite{JSW} for a negative result.\\
On the other hand in \cite{MiI14b} the equivalence of integral ISS and uniform exponential stability is shown for a natural infinite-dimensional version of \eqref{intro:eq3},
 with $A$ generating a $C_{0}$-semigroup and $B:X\times U\to X$ satisfying a Lip\-schitz condition on bounded subsets of $X$ uniformly in the second variable and being bounded in the sense that  $\|B(x,u)\| \lesssim \|x\|\gamma(\|u\|)$ for some $\mathcal{K}$-function $\gamma$ and all $x$ and $u$. As indicated above, the property whether a system is ISS or integral ISS is more subtle when boundary controls are considered and consequently, the involved input operators become unbounded. This also applies for bilinear systems which --- in the presence of boundary control --- cannot be treated as in the references mentioned above. 
 
 

In this article we establish the abstract theory to overcome such issues. 
More precisely, 
we study infinite-dimensional bilinear  control systems of the abstract form
\begin{equation}\label{intro:bilin}\tag{$\Sigma(A,[B_1,B_2],F)$}
\dot{x}(t) =Ax(t)+B_{1}F(x(t),u_{1}(t)) + B_2u_2(t), \ t\geq0,
\end{equation}
where
$A$ generates a $C_0$-semigroup on a Banach space $X$ and $B_{1}$ and $B_{2}$ are possibly unbounded linear operators defined on Banach spaces $\widebar{X}$ and $U_{2}$, respectively. The nonlinearity $F:X\times U_{1}\to \widebar{X}$ is assumed to satisfy a Lipschitz condition and to be bounded. 
In Section \ref{sec2} we present the details of this abstract framework and derive the main result, which, in terms of integral ISS, see also Definition \ref{def:ISS}, reads as follows. 

\begin{mainresult}[Theorem \ref{thm2}]
The bilinear system \ref{eqn1} is integral ISS, if the linear systems $\Sigma(A,[0,B_1],F)$ and $\Sigma(A,[0,B_2],F)$ are integral ISS.
\end{mainresult}
In order to prove this statement we show existence of global mild solutions to \ref{intro:bilin} by classical fixed point arguments under the weak conditions on the operators $B_{1}$, $B_{2}$. 
	 We conclude by applying our abstract result to the example of a bilinearly controlled Fokker--Planck equation with reflective boundary conditions, which has recently appeared in  \cite{BreitenKunisch,Hartmannetal}.

\section{Input-to-state stability for bilinear systems}\label{sec2}
\subsection{System class and notions}\label{sec2.1}

In the following 
we study bilinear  control systems of the form
\begin{equation}\label{eqn1}\tag{$\Sigma(A,[B_1,B_2],F)$}
\arraycolsep=1.3pt\def\arraystretch{1.2}
\begin{array}{ll}\dot{x}(t) &=Ax(t)+B_{1}F(x(t),u_{1}(t)) + B_2u_2(t), \ t\geq0,\\
 x(0)&=x_0,
 \end{array}\hspace{-0.3cm}
\end{equation}
where
\begin{itemize}
\item[$\bullet$] $X$, $\widebar{X}$ and $U_{1}, U_{2}$ are Banach spaces and $x_{0}\in X$,
\item[$\bullet$] $A$ generates a $C_0$-semigroup $(T(t))_{t\ge 0}$ on $X$, 
\item[$\bullet$] the input functions $u_1$ and $u_2$ are locally integrable function with values in ${U_{1}}$ and $U_2$, respectively, that is,  $u_1 \in L^1_{\rm loc}(0,\infty;U_{1})$ and 
$u_2\in L^1_{\rm loc}(0,\infty;U_2)$,
\item[$\bullet$] the operators $B_1$ and $B_2$ are defined on $\widebar{X}$ and $U_2$ respectively. Both operators map
into a space (see below) in which $X$ is densely embedded,
\item[$\bullet$] the nonlinear operator $F:X\times U_1 \to \widebar{X}$ is bounded in the sense that there exists a constant $m>0$ such that 
\begin{equation}\label{eq:bddF}
	\|F(x,u)\|_{\widebar{X}}\leq m \|x\|_X \|u\|_{U_{1}} \qquad \forall x\in X, u\in U_1.
\end{equation}
and Lipschitz continuous in the first variable on bounded subsets of $X$, where the Lipschitz constant
depends on the $U_{1}$-norm of the second argument,
 that is, for all bounded subsets $X_{\text{b}}\subset X$ there exists a constant $L_{X_{\text{b}}}>0$, such that 
\begin{equation}\label{eq:LipschitzF}
\|F(x,u)-F(\widebar{X},u) \|_{\widebar{X}} \leq L_{X_{\text{b}}} \|u\|_{U_1} \|x-\widebar{X}\|_X \qquad \forall x, \tilde{x }\in X_{\text{b}}, u\in U_1,
\end{equation}
\item[$\bullet$] $s\mapsto F(f(s),g(s))$ is measurable for any interval $I$ and measurable functions $f:I\to X$, $g:I\to U_1$,
\item[$\bullet$] we write $\Sigma(A,[0,B_2])$ if $B_{1}=0$ and thus System $\Sigma(A,[0,B_2])$ is linear.
\end{itemize}
Before explaining the details on the assumptions on $B_{1}$ and $B_{2}$ below, we list some examples for functions $F$ and operators that fit our setting.
\begin{enumerate}[label=(\alph*)]
	\item $\widebar{X}=X$, $U=\C$ and $F(x,u)=xu$, 
	\item $\widebar{X}=U=X$, $f\in X^{*}$, $F(x,u)=f(x)u$, 
	\item $\widebar{X}=\C$, $U=X^{*}$, $F(x,u)=\langle x,u\rangle$.
\end{enumerate}


Let $X_{-1}$ be  the completion of $X$ with respect to the norm 
$ \|x\|_{X_{-1}}= \|(\beta -A)^{-1}x\|_X$ 
for some $\beta $ in the resolvent set $\rho(A)$ of $A$. 
For a reflexive Banach space, $X_{-1}$ can be identified with $(D(A^*))'$, the continuous dual of $D(A^*)$ with respect to the  pivot space $X$. 
The operators $B_{1}$ and $B_{2}$ are assumed to map to $X_{-1}$, more precisely, $B_{1}\in  L(\widebar{X},X_{-1})$ and $B_{2}\in  L(U_2,X_{-1})$, where $L(X,Y)$ refers to the bounded linear operators from $X$ to $Y$. Only in the special case that $B_{1}$ or $B_{2}$ are in $L(\widebar{X},X)$ or $L(U_{2},X)$, we say that the respective operator is bounded.
The $C_{0}$-semigroup  $(T(t))_{t\ge 0}$ extends uniquely to a $C_{0}$-semigroup  $(T_{-1}(t))_{t\ge 0}$ on $X_{-1}$ whose generator $A_{-1}$ is the unique extension of $A$ to an operator in $L(X,X_{-1})$, see e.g.\ \cite{engelnagel}. Note that $X_{-1}$ can be viewed as taking the role of a Sobolev space with negative index.
 With the above considerations we may consider System \ref{eqn1} on the Banach space $X_{-1}$.
We want to emphasize that our interest is primarily in the situation where $B_1$ and $B_2$ are not bounded --- something that typically happens if the control enters through point boundary actuation. 
 Note, however, that the assumptions imply that 
``the unboundedness of $B_1$ and $B_2$ is not worse than the one of $A$'' --- which particularly means that if $A\in L(X,X)$ then $B_1\in L(\widebar{X},X)$ and $B_2\in L(U_{2},X)$.  \\
For zero-inputs $u_1$ and $u_2$, the solution theory for System \ref{eqn1} is fully characterized by the property that $A$ generates a $C_{0}$-semigroup as this reduces to solving a linear, homogeneous equation. For non-trivial inputs, the solution concept is a bit more delicate.

\smallskip

More precisely, for given $t_{0},t_{1}\in [0,\infty)$, $t_{0}<t_{1}$, $x_{0}\in X$, $u_1\in L_{loc}^{1}(0,\infty;U_1)$ and $u_{2}\in L_{loc}^{1}(0,\infty;U_{2})$, a continuous function $x:[t_{0},t_{1}] \rightarrow X$ is called a {\em mild solution of} \ref{eqn1}
 {\em on $[t_{0},t_{1}]$} if for all $t\in [t_{0},t_{1}]$,
\begin{equation}
\label{eq:mild}
 x(t) = T(t-t_{0})x_0 + \int_{t_{0}}^t T_{-1}(t-s)[B_1 F(x(s),u_1(s)) + B_2u_2(s)] \, \mathrm{d}s.
\end{equation}
We say that $x:[0,\infty)\to X$ is a {\em global mild solution} or a {\em mild solution on $[0,\infty)$} of \ref{eqn1} if $x|_{[0,t_{1}]}$ is a mild solution on $[0,t_{1}]$ for every $t_{1}>0$. 
 We stress that existence of a mild solution is non-trivial, even when $u_{1}=0$. 
In this case, it is easy to see that $x\in C([0,\infty); X_{-1})$, but not necessarily $x(t)\in X$, $t>0$,  without further assumptions on $B_{2}$.
The existence of a mild solutions to the linear System $\Sigma(A,[0,B_2])$ is  closely related to the notion {\em admissibility} of the operator $B_2$ for the semigroup $(T(t))_{t\ge 0}$ and various sufficient and necessary conditions are available, see e.g.\ Proposition \ref{proplin} and \cite{JNPS16}.
\smallskip

We need the following well-known function classes from Lyapunov theory.
\begin{align*} 
\arraycolsep=1.3pt\def\arraystretch{1.2}
\calK ={}& \{\mu\in C(\Rpn,\Rpn) \:|\: \mu(0)=0, \mu \text{ strictly increasing}\},\\
 \calK_\infty ={}& \{\theta\in\calK \:|\:  \lim_{x\to\infty} \theta(x)=\infty\},\\
\mathcal{L}={}&\{\gamma\in C(\Rpn,\Rpn)\:|\:\gamma \text{ strictly decreasing,}  \lim_{t\to\infty}\gamma(t)=0 \},\\
\mathcal{KL} = {}&\{ \beta:(\Rpn)^{2}\rightarrow\Rpn\: | \: \beta(\cdot,t)\in\calK\ \forall t \geq 0,\beta(s,\cdot)\in\mathcal{L}\ \forall s>0\}.
\end{align*}

The following concept is central in this work. It originates from works by Sontag \cite{Sontag89ISS,Son98}. We refer e.g.\ to \cite{MiroPrie20,MiroWirt18a} for the infinite-dimensional setting.
\begin{definition}\label{def:ISS}
The system \ref{eqn1} is called
\begin{enumerate}[label=(\roman*)]
\item {\em input-to-state stable (ISS)}, if there exist $\beta\in \mathcal{KL}$, $ \mu_1, \mu_2 \in {\calK}_\infty$ such that for every $x_0\in X$, $u_1 \in L^\infty(0,\infty;U_1)$ and $u_2 \in L^\infty(0,\infty;U_2)$ there exists a unique global mild solution $x$ of \ref{eqn1} and for every $t\geq 0$ 
\begin{equation*}
\left\| x(t)\right\| \le \beta(\|x_0\|,t)+  \mu_1 (\|u_1\|_{L^\infty(0,t;U_1)}) + \mu_2 (\|u_2\|_{L^\infty(0,t;U_2)});
\end{equation*}
\item {\em integral input-to-state stable (integral ISS)}, if there exist $\beta\in \mathcal{KL}$, $ \theta_1, \theta_2 \in {\calK}_\infty$ and $\mu_1, \mu_2 \in {\calK}$ such that for every $x_0\in X$, $u_1\in L^\infty(0,\infty;U_1)$ and $u_2 \in L^\infty(0,\infty;U_2)$ System \ref{eqn1} has a unique global mild solution $x$ and for every $t\ge 0$ 
\begin{equation*}
\left\| x(t)\right\| \le \beta(\|x_0\|,t) + \theta_1 \left(\int_0^t \mu_1 (\|u_1(s)\|) \, \mathrm{d}s\right) + \theta_2 \left(\int_0^t \mu_2 (\|u_2(s)\|) \, \mathrm{d}s\right).
\end{equation*}
\end{enumerate}
One may define some mixed type of these definitions like (ISS,integral ISS) (and (integral ISS,ISS)), in the sense that one has an ISS-estimate for $u_1$ and some integral ISS-estimate for $u_2$ (and vice versa).
\end{definition}


\smallskip
Although the terms involving $u_{1}$ and $u_{2}$ on the right-hand-side of the integral ISS estimate do not define norms in general. However, there are function spaces which are naturally linked to integral ISS \cite{JNPS16}.
 In this context we briefly introduce the Orlicz space $E_\Phi(I;Y)$ associated to a so-called Young function $\Phi$ for an interval $I \subset \R$ and a Banach space $Y$ in the \nameref{appendix}. Note that the Orlicz space $E_\Phi$ corresponding to the Young function $\Phi(t)=t^p$, $1<p<\infty$, is isomorphic to $L^p$.
\begin{definition}Let $(T(t))_{t\ge 0}$ be a $C_0$-semigroup.
\begin{enumerate}[label=(\roman*)]
\item We say that $(T(t))_{t\ge 0}$ is {\em of type $(M,\omega)$} if $M\geq1$ and $\omega\in\R$ are such that
\begin{equation} \label{eqn10} \|T(t)\| \le M\e^{-\omega t},\qquad t\ge 0.\end{equation}
\item We say that $(T(t))_{t\ge 0}$ is {\em (uniformly) exponentially stable} if $(T(t))_{t\ge 0}$ is of type $(M,\omega)$ for some $\omega>0$.
\item Let $Z=E_\Phi$ or $Z=L^\infty$. An operator $B \in L(U,X_{-1})$ is called {\em $Z$-admissible} for $(T(t))_{t \ge 0}$,  
if 
 for every $t>0$ and $u\in Z(0,t;U)$ it holds that
  \[\int_0^t T_{-1}(t-s)B u(s) \, \mathrm{d}s\in X.\]

 We will neglect the reference to $(T(t))_{t \ge 0}$  if this is clear from the context.
\end{enumerate}
 \end{definition}
Recall that every $C_0$-semigroup is of type $(M,\omega)$ for some $M\geq1$ and $\omega \in \mathbb{R}$.
Note that any bounded operator $B$ is $Z$-admissible for all $Z$ considered above.
\begin{remark}
\label{rem1}
Let $B\in L(U,X_{-1})$ be $Z$-admissible for $(T(t))_{t\ge 0}$ with $Z=E_\Phi$ or $Z=L^\infty$. Then for any $t>0$ there exists a minimal constant $C_{t,B}>0$ such that 
 \begin{equation}\label{admissibility}
\left\| \int_0^t T_{-1}(t-s)B u(s) \, \mathrm{d}s\right\| \le C_{t,B} \|u\|_{Z(0,t;U)}, \quad u\in Z(0,t;U).
\end{equation}
This is a consequence of the closed graph theorem. Also note that $B$ is $Z$-admissible for $(\e^{\delta t}T(t))_{t\ge 0}$ for any $\delta\in\R$. Furthermore, the function $t\mapsto C_{t,B}$ is non-decreasing and, if $(T(t))_{t\ge 0}$ is exponentially stable, even bounded, that is, $C_{B}:=\sup_{t\ge 0} C_{t,B}<\infty$.

\end{remark}

The following result clarifies the relation between admissibility and (integral) ISS. The interest to study admissibility with respect to Orlicz spaces follows by the natural connection to integral ISS for linear systems, see Proposition \ref{prop:linear} $(iii)$.\\
Note in particular that the existence of mild solutions for $E_{\Phi}$-admissible operators $B_{2}$ is based on the absolute continuity of the Orlicz norm with respect to the length of the interval and the strong continuity of the shift-semigroup on $E_\Phi(I;Y)$ for any interval $I$ and any Banach space $Y$. The latter can be proven by similar methods one uses to prove the strong continuity of the shift-semigroup on $L^p(I;Y)$. 

\begin{proposition}[Prop.~2.10 $\&$ Thm.~3.1 in {\cite{JNPS16}}]\label{prop:linear}
Let $A$ generate the $C_{0}$-semigroup $(T(t))_{t \ge 0}$ on $X$ and $B_{2}\in L(U_{2},X_{-1})$.
\label{proplin}
\begin{enumerate}[label=(\roman*)]
\item If $B_{2}$ is $E_\Phi$-admissible, then for every $x_0 \in X$ and $u_{2} \in E_{\Phi,{\rm loc}}(0,\infty;U_2)$ there exists a unique global mild solution $x$ of System $\Sigma(A,[0,B_2])$, which is given by \eqref{eq:mild} with $B_1=0$.

\item System $\Sigma(A,[0,B_{2}])$ is ISS if and only if $(T(t))_{t\ge 0}$ is  exponentially stable and $B_{2}$ is $L^\infty$-admissible.
\item $\Sigma(A,[0,B_{2}])$ is integral ISS if and only if  $(T(t))_{t \ge 0}$ is exponentially stable and $B_{2}$ is $E_\Phi$-admissible for some Young function $\Phi$.
\end{enumerate}
\end{proposition}


\subsection{Main results}
Whether ISS implies integral ISS for a system $\Sigma(A,[0,B])$  is still an open question. This is true for $B$ bounded, see e.g. \cite[Prop.~2.14]{JNPS16} or \cite{MiI14b}. However, various conditions on $A$ and the input spaces $U$ are available under which integral ISS and ISS are equivalent \cite{JacoSchwZwar19} in the case of boundary control. 

The following proposition proves an estimate between  Orlicz-norms and integral ISS estimates, which will be useful for the proof of the main result. 

\begin{proposition}\label{prop:Orliczintegral ISS}
Let $\Phi$ be a Young function. Then there exist $\mathcal{K}_{\infty}$-functions $\theta$ and $\mu$ such that for any Banach space $U$ and $t>0$,
\begin{equation}\label{eq:Orliczintegral ISSest}
 \|u\|_{E_{\Phi}(0,t;U)} \leq \theta\left(\int_{0}^{t}\mu(\|u(s)\|_{U}) \, \mathrm{d}s\right),\qquad \forall u\in L^\infty(0,t;U).
\end{equation}
Moreover, $\theta$ and $\mu$ can be chosen as
\begin{equation}\label{eq:mu}
		 \mu(x) = \begin{cases}
    \int_{0}^{x}\phi(\sqrt{s}) \, \mathrm{d}s, &x< 1, \\\\
    \frac{ \int_{0}^{1}\phi(\sqrt{s}) \, \mathrm{d}s}{\Phi(1)} \Phi(x^{2}), &x\geq 1,
  \end{cases} 
 \end{equation}
where $\phi$ equals the right-derivative of $\Phi$ a.e.\ and, for $\alpha >0$,
\begin{equation*}
	\theta (\alpha) = \sup\left\{ \|u\|_{E_{\Phi}(0,t;U)} \Bigm| u \in L^\infty(0,t;U),\,t \geq 0,\,\int_0^t \mu(\|u(s)\|_U) \,\mathrm{d}s \leq \alpha \right\},
\end{equation*}
with $\theta(0)=0$. \\
If $\Phi$ satisfies the $\Delta_2$-condition (c.f.~\nameref{appendix}), then $\mu=\Phi$ can be chosen as well.
\end{proposition}

\begin{proof}
Note that we only need to show that $\mu$ and $\theta$ define $\mathcal{K}_{\infty}$-functions since \eqref{eq:Orliczintegral ISSest} is immediate from the definition of $\theta$.
The proof is similar in spirit to an argument used in \cite[Proof of Thm.~1]{NabSchwe}, with the crucial fact being that $\mu$ defined by \eqref{eq:mu} defines a Young function such that
 \[\Phi\leq\mu\quad\text{and}\quad \sup_{x>0}\frac{\Phi(cx)}{\mu(x)}<\infty,\] 
for all $c>0$, see \cite[Lem.~1]{NabSchwe}. In the special case that $\Phi$ satisfies the $\Delta_{2}$-condition (with $s_{0}=0$), the above properties also hold for $\mu=\Phi$, by the defining properties of the $\Delta_{2}$-condition. 
This implies that whenever a sequence $(f_{n})_{n\in\mathbb{N}}$ with  $f_n\in L^\infty(0,t_n;U)$ is such that 
\[ \lim_{n\to\infty}\int_{0}^{t_n}\mu(\|f_{n}(s)\|_{U})\, \mathrm{d}s=0,\]
it follows that $\lim_{n\to\infty}\|f_{n}\|_{E_{\Phi}(0,t_n;U)}=0$, see \cite[Lem.~2]{NabSchwe}. Clearly, $\mu$ is a $\mathcal{K}_{\infty}$-function, since $\mu$ is a Young function. Therefore, it remains to consider $\theta$. It is easy to see that $\theta$ is well-defined, non-decreasing and unbounded, whence we are left to show continuity. Moreover, since $\theta(\alpha)$ is of the form $\sup M_{\alpha}$ with nested sets $(M_{\alpha})_{\alpha>0}$, it follows that $\theta$ is right-continuous on $(0,\infty)$. To see that $\theta$ is continuous at $\alpha=0$, let $(\alpha_{n})_{n}$ be a decreasing sequence of positive numbers with $\lim_{n\to\infty}\alpha_{n}=0$ and for every $n \in \mathbb{N}$ let $u_{n}\in L^\infty(0,t_n;U)$ be such that $\int_{0}^{t_n}\mu(\|u_{n}(s)\|)\, \mathrm{d}s\leq \alpha_{n}$ and $0\leq \theta(\alpha_{n})-\|u_{n}\|_{E_{\Phi}(0,t_n;U)}<\frac{1}{n}$. By the above mentioned argument, we can conclude that $\|u_{n}\|_{E_{\Phi}(0,t_n;U)}$ converges to $0$ as $n\to\infty$. Thus, $\lim_{n\to\infty}\theta(\alpha_{n})=0$.
\newline We finish the proof by showing that $\theta$ is left-continuous on $(0,\infty)$. Now let $\alpha>0$, $\alpha_{n}\nearrow\alpha$ and let $u_{k}\in L^\infty(0,t_k;U)$, $k\in\mathbb N$, such that 
\[\int_{0}^{t_k}\mu(\|u_{k}(s)\|) \, \mathrm{d}s\leq \alpha\quad \text{and}\quad \lim_{k\to\infty}\theta(\alpha)-\|u_{k}\|_{E_{\Phi}(0,t_k;U)}=0.\] For every $n\in\mathbb{N}$, we aim to find $\tilde{u}_{n}\in L^\infty(0,t_n;U)$ such that  $\int_{0}^{t_n}\mu(\|\tilde{u}_{n}(s)\|) \, \mathrm{d}s\leq \alpha_{n}$ and $\lim_{n\to\infty}\|u_{n}-\tilde{u}_{n}\|_{E_{\Phi}(0,t_n;U)}=0$. Indeed, then 
	\begin{align*}
		\theta(\alpha)-\theta(\alpha_{n})\leq{}&\theta(\alpha)-\|\tilde{u}_{n}\|_{E_{\Phi}(0,t_n;U)}\\
			\leq{}&\theta(\alpha)-\|{u}_{n}\|_{E_{\Phi}(0,t_n;U)}+\|u_{n}-\tilde{u}_{n}\|_{E_{\Phi}(0,t_n;U)}
	\end{align*}
tends to $0$ as $n\to\infty$, which shows left-continuity. 
We define $\tilde{u}_{n}:=u_{n}\chi_{M_{n}}$
where the measurable set $M_{n}$ is chosen such that 
\[\int_{M_{n}}\mu(\|u_{n}(s)\|) \, \mathrm{d}s=\alpha_{n}, \quad \text{if }\int_{0}^{t_n}\mu(\|u_{n}(s)\|) \, \mathrm{d}s\geq\alpha_{n},\]
or $M_{n}=(0,t_n)$ otherwise. 
 It follows that 
\begin{align*}
\int_{0}^{t_n}\mu(\|u_{n}(s)-\tilde{u}_{n}(s)\|_{U})\, \mathrm{d}s={}&\int_{0}^{t_n}\mu(\|u_{n}(s)\|_{U}) \, \mathrm{d}s-\int_{M_{n}}\mu(\|{u}_{n}(s)\|_{U})\, \mathrm{d}s\\
\leq{}&\alpha-\alpha_{n}.
\end{align*}
Thus, using the argument from the beginning of the proof again, we infer that $\|u_{n}-\tilde{u}_{n}\|_{E_{\Phi}(0,t_n;U)}\to0$ as $n\to\infty$. This concludes the proof.
\end{proof}

Combining Proposition \ref{prop:Orliczintegral ISS} with \cite[Prop~2.10]{JNPS16} allows us to fomulate the following result:
\begin{corollary}
If System $\Sigma(A,[0,B_2])$ possesses a unique mild solution $x$ for every $x_0 \in X$ and $u_2 \in  L^\infty(0,\infty;U_2) $, then the following statements are equivalent.
\begin{enumerate}[label=(\roman*)]
	\item There exist functions $\beta \in \mathcal{KL}$ and $\mu_2 \in \mathcal{K}_\infty$ such that
\begin{equation}
\lVert x(t) \rVert \leq \beta(\lVert x_0 \rVert,t) + \mu_2(\lVert u_2 \rVert_{E_\Phi(0,t;U_2)})
\label{eq:OrliczISS}
\end{equation}
holds for all $t \geq 0$ and $u_2 \in  L^\infty(0,\infty;U_2)$.
\item There exist functions $\beta \in \mathcal{KL}$, $\theta_2 \in \mathcal{K}$ and $\mu_2 \in \mathcal{K}_\infty$ such that
\begin{equation}
\lVert x(t) \rVert \leq \beta( \lVert x_0 \rVert , t ) + \theta_2 \left( \int_0^t \mu_2(\lVert u_2(s) \rVert_{U_2}) \, \mathrm{d} s \right)
\label{eq:Orliczintegral ISS}
\end{equation}
holds for all $t \geq 0$ and $u_2 \in L^\infty(0,\infty;U_2)$.
\end{enumerate}
\end{corollary}

\begin{remark}\label{constructionofmuandtheta}
Let us make the following comments on the construction of $\mu$ and $\theta$ in Proposition \ref{prop:Orliczintegral ISS}.
\begin{enumerate}
	\item If $\Phi(s)=s^{p}$, $s>0$, then $\mu(s)=s^{p}$ and it is not hard to see that, up to a constant, $\theta(r)$ is given by $\Phi^{-1}(r)=r^{\frac{1}{p}}$. This shows that the choice of $\theta$ is rather natural.
	\item 
	With similar techniques as in the proof of Proposition \ref{prop:Orliczintegral ISS}, it has been shown in \cite{JNPS16,NabSchwe} that if a linear system $\Sigma(A,[0,B])$ satisfies  \eqref{eq:OrliczISS}, then it is integral ISS with the estimate
	\[ \|x(t)\|\leq \beta(\|x_{0}\|,t)+\theta\left(\int_{0}^{\infty}\mu(\|u(s)\|_{U}) \, \mathrm{d}s\right),\]
	where 
	\begin{align}
	 \notag\theta (\alpha) = {}\sup&\left\{ \left\|\int_{0}^{t}T_{-1}(s)Bu(s) \, \mathrm{d}s\right\| \Bigm|\right.\\{}\quad &\left. 
	  u \in L^\infty (0,t;U),\,t \geq 0,\,\int_0^t \mu(\|u(s)\|_U) \, \mathrm{d}s \leq \alpha \right\}.\label{eq:thetadiff}
	 \end{align}
	 Proposition \ref{prop:Orliczintegral ISS} shows that $\theta$ can actually be chosen independent of the semigroup $(T(t))_{t \geq 0 }$ and $B$ provided the system is $E_{\Phi}$-admissible (which, however, depends on $(T(t))_{t \geq 0}$ and $B$, of course). In some sense, this fact simplifies the proofs in \cite{JNPS16,NabSchwe}. On the other hand, the choice of $\theta$ based on \eqref{eq:thetadiff} is more refined; in case the system was even $E_{\Psi}$-admissible with some $\Psi\leq\Phi$, this would affect the choice of $\theta$, even if $\mu$ is constructed from $\Phi$ only.
\end{enumerate}
	\end{remark}
	
%

In contrast to linear systems, the existence of mild solutions is less clear for bilinear control systems of the form \ref{eqn1}. 
 
Sontag \cite{Son98} showed that  finite-dimensional  bilinear systems are hardly ever ISS, but integral ISS if and only if the semigroup is exponentially stable. 
 In \cite{MiI14b} it was shown that this result generalizes to infinite-dimensional bilinear systems
provided that $B_{1}$ and $B_{2}$ are bounded operators and $\widebar{X}=X$. The following results give sufficient conditions for integral ISS and some combination of ISS and integral ISS of \ref{eqn1}.
We start with a result on existence of local solutions to \ref{eqn1}. The proof involves typical arguments in the context of mild solutions for semilinear equations.

A similar result for the existence of the unique 
 mild solution as in the following Lemma \ref{thm1}  were proved under slightly stronger conditions in \cite{Berra09} for $L^{p}$-admissible $B_{1}$, scalar-valued inputs $u_1$, $F(x,u_1)=u_1 x$ and $B_2=0$. Our condition is more natural as the same condition guarantees the existence of continuous (and unique) global mild solutions of the linear systems $\Sigma(A,[0,B_1])$ and $\Sigma(A,[0,B_2])$, see Proposition \ref{proplin}.

\begin{lemma}\label{thm1}
Let $A$ generate a $C_{0}$-semigroup $(T(t))_{t \ge 0}$ on $X$. Suppose that  $B_{1}\in L(\widebar{X},X_{-1})$ is $E_\Phi$-admissible and that $B_{2}\in L(U_{2},X_{-1})$ is $E_\Psi$-admissible.
Then  for every $t_0\geq 0$, $x_0\in X$, $u_1\in E_\Phi(0,\infty;U_1)$ and $u_2\in E_\Psi(0,\infty;U_2)$ there exists $t_1>t_0$ such that System \ref{eqn1} possesses a unique mild solution $x$ on $[t_0,t_1]$. \\Moreover, if $t_{\mathrm{max}} > t_0$ denotes the supremum of all $t_1>t_0$ such that System \ref{eqn1} has a unique mild solution $x$ on $[t_0,t_1]$, then $t_{\mathrm{max}} < \infty$ implies that
\[ \lim_{t \to t_{\mathrm{max}}} \|x(t)\| =  \infty.\]
\end{lemma}

\begin{proof}
We first show that for every $t_0 \geq 0, x_0\in X$, $u_1\in  E_\Phi(0,\infty;U_1) $ and $u_2\in  E_\Psi(0,\infty;U_2)$ there exists $t_1>t_{0}$ such that System \ref{eqn1} possesses a unique mild solution on $[t_0,t_1]$ with initial condition $x_0$ and input functions $u_1$ and $u_2$. Moreover, we show that $t_{1}=t_{0}+\delta$ can be chosen such that $\delta$ is independent for any bounded sets of initial data $x_{0}$ and  $t_{0}$.  Let $T>0$, $r>0$, $u_1\in  E_\Phi(0,\infty;U_1) $ and $u_2\in  E_\Psi(0,\infty;U_2)$ be arbitrarily. We first recall the following property of Orlicz spaces. For any $\varepsilon>0$ there exists $\delta>0$ such that 
\begin{equation}\label{:eq:absolutenorm}
\max\{\|u_1\|_{E_\Phi(t,t+\delta;U_1)}, \|u_2\|_{E_\Psi(t,t+\delta;U_2)}\}< \varepsilon,\qquad \forall t \ge 0,
\end{equation}
see e.g. \cite[Thm.~3.15.6]{Kufner}.
Let $t_{0}\in[0,T]$, $t_1>t_0$ and $x_{0}\in K_{r}(0)=\{x\in X\colon \|x\|\leq r\}$ and define the mapping 
\begin{align*}
\Phi_{t_0,t_1}&:   C([t_0,t_1];X)\rightarrow C([t_0,t_1];X)\\
( \Phi_{t_0,t_1}(x))(t)&:=T(t-t_0)x_0+\int_{t_0}^t T_{-1}(t-s)[B_1F(x(s),u_1(s))+ B_2u_2(s)] \, \mathrm{d}s.
\end{align*}
The strong continuity of $(T(t))_{t\ge 0}$ and Proposition \ref{proplin} imply that $\Phi_{t_0,t_1}$ is well-defined, that is,  $\Phi_{t_0,t_1}(x)\in C([t_0,t_1];X)$ for every  $ x\in  C([t_0,t_1];X)$.
Note that we applied Proposition \ref{proplin} twice: to System $\Sigma(A,[0,B_2])$ with input $u_2$ and to System $\Sigma(A,[0,B_1])$ with input  $F(x(\cdot),u_1(\cdot))$, where we set $u_1,u_2,x$ zero on $(0,t_0)$.

Let $M \geq 1$ and $\omega \in \R$ be such that $\|T(t)\|\leq M \e^{-\omega t}$ for all $t \geq 0$ and choose $k=4Mr + 2M$.  Set
\[M_k(t_0,t_1):= \{ x \in C([t_0,t_1];X) \mid \|x\|_{C([t_0,t_1];X)} \leq k \}.\]
We will show next that $t_1$ can be chosen  such that $\Phi_{t_0,t_1}$ maps  $M_k(t_0,t_1)$ to $M_k(t_0,t_1)$ and is contractive on this set. 
 Let $C_{t,B_{1}}$ and $C_{t,B_{2}}$ refer to the admissibility constants such that \eqref{admissibility} holds for $B_1$ and $B_2$ which can be chosen non-decreasing in $t$. Let $m$ be the boundedness constant of $F$ from \eqref{eq:bddF} and let $L_{K_{k}(0)}$ be the Lipschitz constant of $F$ such that \eqref{eq:LipschitzF} holds for the bounded set $X_{b}=\{ x(t) \mid x \in M_k(t_0,t_1), t \in [t_0,t_1] \} \subset X$ which is equal to $K_{k}(0)=\{x\in X\colon \|x\|\leq k\}$. 
Now, let $t_1 = t_0 + \delta$ with $\delta \in (0,1)$ be chosen such that for all $t_{0}\in[0,T]$,
\begin{enumerate}[label=(\roman*)]
	\item $\e^{\omega (t_1-t_0)} =\e^{\omega \delta} \leq 2 $, 
	\item $m C_{T+1,B_{1}}\|u_1\|_{E_\Phi(t_{0},t_{0}+\delta;U_1)} \leq \frac{1}{2}$,
	\item $C_{T,B_{2}} \|u_2\|_{E_\Psi(t_{0},t_{0}+\delta;U_2)}\leq M$ and
	\item $C_{T+\delta,B_{1}} L_{K_k(0)} \|u_1\|_{E_\Phi(t_{0},t_{0}+\delta;U_1)} < 1$
\end{enumerate}
holds, where we used \eqref{:eq:absolutenorm} in (ii)-(iv). Note that apart from the parameters of the operators $B_{1},B_{2}, A,F$, the choice of $\delta$ only depends on $r$ and $T$,  where the $r$-dependence of $\delta$ arises from the $r$-dependence of $k$. 
 It follows that for all $t_{0}\in[0,T]$, $x\in M_k(t_0,t_1)$ and $x_0\in K_r(0)$
\begin{align*}
&\|\Phi_{t_0,t_1}(x)\|_{C([t_0,t_1];X)}\\
&\leq M \e^{\omega (t_1-t_0)} \|x_0\| + C_{t_{1},B_{1}} \| F(x,u_1) \|_{E_\Phi(t_0,t_1;\widebar{X})} + C_{t_{1},B_{2}} \|u_2\|_{E_\Psi(t_0,t_1;U_2)} \\
&\leq 2M \|x_0\| + m C_{t_1,B_{1}} \|u_1\|_{E_\Phi(t_0,t_1;U_1)} \|x\|_{C([t_0,t_1];X)} + M \\
&\leq k,
\end{align*}
where we used admissibility in the first inequality and \eqref{eq:bddF} in the second inequality as well as the monotonicity of the Orlicz norm in both estimates. Hence, $\Phi_{t_0,t_1}$ maps $M_k(t_0,t_1)$ to $M_k(t_0,t_1)$.
The contractivity follows by $\rm{(iv)}$ since
\begin{eqnarray*}
\lefteqn{\|\Phi_{t_0,t_1}(x)-\Phi_{t_0,t_1}(\tilde x)\|_{C([t_0,t_1];X)}}\nonumber\\
&\le & \sup_{t\in[t_0,t_1]} \left\|\int_{t_0}^{t}T(t-s)B_1[F(x(s),u_1(s)) -F(\tilde{x}(s),u_1(s))]\, \mathrm{d}s\right\| \label{eqn5}\\
&\le & C_{t_{1},B_{1}} L_{K_k(0)} \|u_1\|_{E_\Phi(t_0,t_1;U_1)}  \|x-\tilde x\|_{C([t_0,t_1];X)},
\end{eqnarray*}
where we used again admissibility, the Lipschitz property of $F$ and the monotonicity of the Orlicz norm. By Banach's fixed-point theorem, we conclude that  System \ref{eqn1} possesses a unique mild solution on $[t_0,t_1]$ with initial condition $x_0$ and input functions $u_1$ and $u_2$.\\
Now let $t_{\mathrm{max}}$ be the supremum of all $t_{1}$ such that there exists a mild solution $x$ of \ref{eqn1} on $[t_0,t_1]$ for every $t_1 < t_{\mathrm{max}}$, where  $x_0 \in X$, $u_1 \in E_\Phi(0,\infty;U_1)$ and $u_2 \in E_\Psi(0,\infty;U_2)$ are given.  Suppose that $t_{\mathrm{max}}$ is finite. We will show, that then $\lim_{t\to t_{\mathrm{max}}}\|x(t)\|=\infty$. If this is not the case, we have
\[ r = \sup_{t \in [t_0,t_{\mathrm{max}}]} \|x(t)\| < \infty.\]

Let $(t_n)_{n \in \mathbb{N}}$ be a sequence of positive real numbers converging to $t_{\mathrm{max}}$ from below. 
Since $t_{n}\in[0,t_{\mathrm{max}}]$ and $\|x(t_{n})\|\leq r$ for all $n\in\mathbb{N}$, there exists $\delta>0$ independent of $n\in\mathbb{N}$ such that the equation
\begin{equation*}
\begin{array}{ll}
\dot{y}(t) = Ay(t) + B_1F(y(t),u_1(t)) + B_2u_2(t), \\
y(t_n) = x(t_n).
\end{array}
\end{equation*}
has a mild solution $y$ on $[t_{n},t_{n}+\delta]$.
Therefore, we can extend $x$ by $x(t) =y(t)$, $t \in (t_n,t_n+ \delta]$, to a solution of \ref{eqn1} on $[t_0,t_n+ \delta]$. This contradicts the maximality of $t_{\mathrm{max}}$ and hence, $x$ has to be unbounded in $t_{\mathrm{max}}$. 
\end{proof}


\begin{theorem}\label{thm2}
Suppose that the linear systems $\Sigma(A,[0,B_1])$ and $\Sigma(A,[0,B_2])$ are integral ISS, then the bilinear system \ref{eqn1} is integral ISS and (integral ISS,ISS).\\
The assumption that $\Sigma(A,[0,B_2])$ is integral ISS is necessary.
%


 
\end{theorem} 

\begin{proof} 
The necessity of $\Sigma(A,[0,B_2])$ being integral ISS follows by setting $u_1=0$ in the bilinear system.\\ 
Proposition \ref{prop:linear} says that integral ISS of the linear systems is equivalent to the exponential stability of the semigroup $(T(t))_{t \geq 0}$ generated by $A$ and the admissibility of the control operators $B_1$ and $B_2$ with respect to some Orlicz spaces $E_\Phi$ and $E_\Psi$, respectively.\\
Using this characterization, we will give the proof in two steps. At first we prove the existence of a global mild solution $x$ of \ref{eqn1} (which does not need the exponential stability of $(T(t))_{t \geq0}$). Afterwards we prove the (integral) ISS properties.\\ 
{\bf STEP I.} Let $(M,\omega)$ denote the type of $(T(t))_{t \geq0}$.
By Remark \ref{rem1} there exist $C_{t,B_1},C_{t,B_2}>0$ such that for every $t\ge 0$, $y\in E_\Phi(0,\infty;\widebar{X})$ and $\tilde{y} \in E_\Psi(0,\infty;U_2)$  we have
\[  \left\|\int_0^t \e^{\frac{\omega}{2}( t-s)} T_{-1}(t-s)B_1y(s) \, \mathrm{d}s\right\| \le C_{t,B_1} \| y \|_{E_\Phi(0,t;\widebar{X})}\]
and
\begin{align*} 
 \left\|\int_0^t \e^{\frac{\omega}{2}( t-s)} T_{-1}(t-s)B_2 \tilde{y}(s)\, \mathrm{d}s\right\| \le{}&  C_{t,B_2} \left\|\tilde{y} \right\|_{E_\Psi(0,t;U_2)}.
 \end{align*}
Let $x_0 \in X$, $u_1 \in E_\Phi(0,\infty;U_1)$ and $u_2 \in E_\Psi(0,\infty;U_2)$ and let $t_{\mathrm{max}}$ be the supremum over all $t_1$ such that \ref{eqn1} possesses a unique $x$ mild solution on $[0,t_1]$. Lemma \ref{thm1} yields $t_{\mathrm{max}}>0$.
For $t \in [0,t_{\mathrm{max}})$ we have that
\pagebreak[1]
\begin{align}
\|&x(t)\| \nonumber\\
&= \left\|T(t)x_0+ \int_0^t T_{-1}(t-s)B_1 F(x(s),u_1(s)) \, \mathrm{d}s+ \int_0^t T_{-1}(t-s)B_2u_2(s) \, \mathrm{d}s \right\| \notag \\
& \le  \left\|T(t)x_0\right\|+\e^{-\frac{\omega}{2} t}\left\|\int_0^t  \e^{\frac{\omega}{2}( t-s)}T_{-1}(t-s)B_1 (\e^{\frac{\omega}{2} s}F(x(s),u_1(s))) \, \mathrm{d}s\right\| \notag \\
& \qquad \qquad \quad \mspace{5mu} +\e^{-\frac{\omega}{2}t}\left\|\int_0^t \e^{\frac{\omega}{2}( t-s)} T_{-1}(t-s)B_2 \e^{\frac{\omega}{2}s}u_{2}(s)\, \mathrm{d}s\right\|\notag\\
&\le  M \e^{-\omega t}\|x_0\| 
+ C_{t,B_1} \e^{-\frac{\omega}{2} t} \| \e^{\frac{\omega}{2} \cdot }F(x(\cdot),u_1(\cdot))\|_{E_\Phi(0,t;\widebar{X})} +C_{\omega,u_{2},t}\label{Z_2-admissibility},
\end{align}
where $C_{\omega,u_{2},t}=C_{t,B_2}\e^{-\frac{\omega}{2}t} \left\| \e^{\frac{\omega}{2}\cdot}u_2 \right\|_{E_\Psi(0,t;U_2)}$.
The $\|\cdot\|_{E_{\Phi}}$-norm  in the second term can be estimated by the boundedness of $F$,
\[ \| \e^{\frac{\omega}{2} \cdot }F(x(\cdot),u_1(\cdot))\|_{E_\Phi(0,t;\widebar{X})} \leq m \left\| \,  \|u_1(\cdot)\| \, \e^{\frac{\omega}{2} \cdot }  \,\|x(\cdot)\| \, \right\|_{E_\Phi(0,t;\C)}.\]
We pass over to the equivalent norm on $E_{\Phi}$ given in the \nameref{appendix}, \eqref{aequinorm}. Therefore, for $\varepsilon>0$ there exists a function $g \in L_{\tilde{\Phi}}(0,t;\C)$ with $\| g\|_{L_{\tilde{\Phi}}(0,t;\C)} \leq 1$ such that 
\[ \left\| \|u_1(\cdot)\| \e^{\frac{\omega}{2} \cdot }\|x(\cdot)\| \right\|_{E_\Phi(0,t;\C)} \leq \int_0^t \|u_1(s)\| \,  |g(s)| \left( \e^{\frac{\omega}{2} s} \|x(s)\| \right) \, \mathrm{d}s + \varepsilon.\]
 Hence, by combining this with \eqref{Z_2-admissibility} gives
\begin{eqnarray*}
\e^{\frac{\omega}{2}  t}\|x(t)\| 
&\le& M \e^{- \frac{\omega}{2}t}\|x_0\|  + m C_{t,B_1} \varepsilon + \e^{\frac{\omega}{2} t}C_{\omega,u_{2},t}\\
&&+~m C_{t,B_1} \int_0^t \| u_1(s)\| \, |g(s)| \left( \e^{\frac{\omega}{2} s} \| x(s) \| \right) \, \mathrm{d}s.
\end{eqnarray*}
Setting $\alpha(t) := M \e^{- \frac{\omega}{2}t}\|x_0\|  + m C_{t,B_1} \varepsilon + \e^{\frac{\omega}{2} t}C_{\omega,u_{2},t}$, Gronwall's inequality implies that
\begin{align*}
\e^{\frac{\omega}{2} t} \|x(t)\|
&\leq  \alpha(t) + m C_{t,B_1} \int_0^t \alpha(s) \| u_1(s)\| \, |g(s)| \e^{\left(m C_{t,B_1} \int_s^t \|u_1(r)\| \, |g(r)| \, \mathrm{d}r \right)}  \mathrm{d}s \\
&\leq\alpha(t) + \left( M \|x_0\| \sup_{r \in [0,t]} \e^{-\frac{\omega}{2}r} + m C_{t,B_1} \varepsilon+ \e^{\frac{\omega}{2} t}C_{\omega,u_{2},t} \right) \\
&\qquad\qquad\quad\cdot 2m C_{t,B_1} \|u_1\|_{E_\Phi(0,t;U_1)} \e^{2m C_{t,B_1} \|u_1\|_{E_\Phi(0,t;U_1)}},
\end{align*}
where we used the generalized H\"older inequality for Orlicz spaces, see \eqref{Holder} in the \nameref{appendix}. 
Thus, by letting $\varepsilon$ tend to $0$, multiplying with $\e^{-\frac{\omega}{2}t}$ and using $ab \leq \frac{1}{2}a^2+\frac{1}{2}b^2$ for $a,b\in\R$, we obtain
\begin{eqnarray*}
 \|x(t)\|
&\le&  M \e^{-\omega t} \|x_0\| + \tfrac{1}{2} M^2 \e^{-\omega t} \sup_{r \in [0,t]} \e^{-\omega r}  \|x_0\|^2  \\
&&{}+4 m^2 C_{t,B_1}^2 \|u_1\|_{E_\Phi(0,t;U_1)}^2 \e^{4 m C_{t,B_1} \|u_1\|_{E_\Phi(0,t;U_1)}}\\
&&{}+ C_{\omega,u_{2},t} + \tfrac{1}{2} C_{\omega,u_{2},t}^2,
 \end{eqnarray*}
by monotonicity of the Orlicz norm,
\begin{equation*}
	 \|\e^{\frac{\omega}{2}\cdot}u_2 \|_{E_\Psi(0,t;U_2)}\le \sup_{r\in [0,t]}\e^{\frac{\omega}{2}r} \| u_2 \|_{E_\Psi(0,t;U_2)}.
\end{equation*}
Thus, we have shown
\begin{align}
\!\!\|x(t)\|\!\label{bilinearISSestimate}
\leq {}&\beta(\|x_{0}\|,t)\!+\!\gamma_{1}(C_{t,B_1}\!\|u_1\|_{E_\Phi(0,t)})\!+\gamma_{2}(C_{t,B_2}\e^{-\frac{\omega}{2}t}\| \e^{\frac{\omega}{2}\cdot}u_2 \|_{E_\Psi(0,t)})\\
\!\leq{}&\beta(\|x_{0}\|,t)\!+\!\gamma_{1}(C_{t,B_1}\!\|u_1\|_{E_\Phi(0,t)})\!+\gamma_{2}(C_{t,B_2} \!\sup_{r \in [0,t]}\e^{-\frac{\omega}{2}r} \|u_2 \|_{E_\Psi(0,t)}),\notag
\end{align}
for all $u_1 \in E_\Phi(0,\infty;U_1)$, $u_2 \in E_\Psi(0,\infty;U_2)$ and functions $\beta \in \mathcal{KL}$ and $\gamma_1, \gamma_2 \in \mathcal{K}_\infty$, which can be choosen as
\begin{align*}
	\beta(s,t)={}&M \e^{-\omega t}s + \tfrac{1}{2} M^2 \e^{-\omega t} s^{2} \sup_{r \in [0,t]}\e^{-\omega r},\\	
	\gamma_{1}(s)={}&4 m^2 s^{2}  \e^{ 4 m s}, \\
	\gamma_{2}(s)={}&s + \tfrac{1}{2} s^2.
\end{align*}
 Moreover, the mild solution exists on $[0,\infty)$. Indeed, if this is not the case, we have $t_{\mathrm{max}}<\infty$ and Lemma \ref{thm1} implies that $x$ is unbounded in $t_{\mathrm{max}}$. This contradicts $\eqref{bilinearISSestimate}$ since the right-hand-side is uniformly bounded in $t$ on finite intervals $[0,t_{\mathrm{max}})$.\\ 
 \rm{\bf STEP II.} Since we are dealing with an exponentially stable semigroup, Remark \ref{rem1} implies that $C_{t,B_1}$ and $C_{t,B_2}$ are uniformly bounded in $t$ and we can choose $\omega>0$. Hence, \eqref{bilinearISSestimate} yields for all $u_1 \in E_\Phi(0,\infty;U_1)$ and $u_2 \in E_\Psi(0,\infty;U_2)$ that
 \begin{align*}
 \|x(t)\| &\leq \beta(\|x_{0}\|,t)+\gamma_{1}\left(C_{B_{1}}~\|u_1\|_{E_\Phi(0,t;U_1)}\right)+\gamma_{2}\left(C_{B_2}~\| u_2 \|_{E_\Psi(0,t;U_2)}\right)
 \end{align*}
with $C_{B_i}= \sup_{t \geq 0}C_{t,B_i}$, $i=1,2$.
 
 Using Proposition \ref{prop:Orliczintegral ISS} for $u_1$ and $u_2$, we have shown that \ref{eqn1} is integral ISS since $L^\infty$ is contained in any Orlicz space on compact intervals. If we apply Proposition \ref{prop:Orliczintegral ISS} only for $u_1$ in \eqref{bilinearISSestimate},  \ref{eqn1} is (integral ISS,ISS) by realizing that there exists a constant $C>0$ such that 
\begin{equation}\label{eq:help}
  \e^{-\frac{\omega}{2}t}\| \e^{\frac{\omega}{2}\cdot}u_2 \|_{E_\Psi(0,t;U_2)} \leq C~\|u_{2}\|_{L^\infty(0,t;U_2)},
  \end{equation}
for all $u_{2}\in L^{\infty}(0,\infty;U_{2})$ and $t>0$. To see this, let $\varepsilon>0$ such that $\Psi(x)\leq x$ for all $x\in(0,\varepsilon)$, which exists by the property that $\lim_{s\to 0}\frac{\Psi(s)}{s}=0$. Therefore, choosing $C=\max\{\frac{1}{\epsilon},\frac{2}{\omega}\}$,
\[\int_{0}^{t}\Psi\left(C^{-1}\e^{-\frac{\omega}{2}s}\right)\mathrm{d}s\leq \int_{0}^{t}C^{-1}\e^{-\frac{\omega}{2}s}\mathrm{d}s\leq \frac{2}{C\omega}\leq1.\]
 This implies that 
\[
\int_{0}^{t}\Psi\left(\frac{\e^{\frac{\omega}{2}s}\|u_2(s)\|}{C\e^{\frac{\omega}{2}t}\|u_2\|_{L^\infty(0,t;U_2)}}\right)\mathrm{d}s\leq\int_{0}^{t}\Psi\left(C^{-1}\e^{\frac{\omega}{2}(s-t)}\right)\mathrm{d}s\le1,
\]
by the definition of the $E_{\Psi}$-norm.
\end{proof}

The assumption in Theorem \ref{thm2} that System $\Sigma(A,[0,B_1])$ is integral ISS is not  necessary as the choice $F=0$ shows.

\begin{remark}\label{rem:aftermainthm}
\begin{enumerate}[label=(\arabic*)]
\item In Theorem \ref{thm2} one cannot expect the bilinear systems to be ISS as the trivial finite-dimensional example $\dot{x}=-x+u_{1}x$ shows.
\item Using the definitions of $\gamma_1$, $\gamma_2$ after \eqref{bilinearISSestimate} and the definitions of $\mu$ and $\theta$ from Proposition \ref{prop:Orliczintegral ISS}, up to constants the functions $\mu_1$, $\mu_2$, $\theta_1$ and $\theta_2$ in the integral ISS estimate for \ref{eqn1} can be given explicitly.
\item The proof of Theorem \ref{thm2} is easier in the case that the Orlicz spaces are $L^p$ spaces, since the $L^p$-norm is already an integral of the form we are seeking for in the integral ISS estimate (c.f.~Definition \ref{def:ISS}).
\item Note that the assumptions of Lemma \ref{thm1} already yield that the unique mild solution is global. This is the first step of the proof of Theorem \ref{thm2}.
\end{enumerate}

\end{remark}

In order to investigate integral ISS, it is thus sufficient to check that the linear systems $\Sigma(A,[0,B_1])$ and $\Sigma(A,[0,B_2])$ are integral ISS, or by Proposition \ref{prop:linear} equivalently, that $A$ generates an exponential stable $C_0$-semigroup and the control operators $B_1$ and $B_2$ are admissible.
Note that 
there are control  operators $B$ which are  $E_\Phi$-admissible for some Young function $\Phi$ but not $L^p$-admissible for any $p \in [1,\infty)$.
In the context of linear systems, such an example was already given in \cite[Ex.~5.2]{JNPS16} for an operator $B$ defined on $\mathbb{C}$ using the connection between a Carleson-measure criterion and admissibility stated in \cite{JNPS16}, see also \cite{JPP14}. 
The next example  extends this result to control operators  defined on $\widebar{X}$. 

\begin{example}
Let $X= \ell^2(\mathbb{N})$ and define
$F:X\times \C\to X$, by $F(x,u):=ux$
and the diagonal operators
\begin{equation*}
Ae_{n}=-2^n e_{n},\quad Be_{n}=\frac{2^n}{n} e_{n},\qquad n\in \mathbb{N},
\end{equation*}
where $(e_n)_{n \mathbb{N}}$ is the canonical basis of $X$ and $A$ is defined on its maximal domain.\\
The general assumptions  of Section \ref{sec2.1} are satisfied with $B_1=B$ and $B_2=0$. Let  $x=(\frac{1}{n})_{n \in \mathbb{N}} \in X$ and $p \in [1, \infty)$. Following \cite[Ex.~5.2]{JNPS16}, the operator $b=Bx$ defined on $\mathbb{C}$ is not $L^p$-admissible.  Hence, $B$ is not $L^p$-admissible.\\
Next, we show that $B$ is $E_\Phi$-admissibile, where $\Phi$ is the complementary Young function to $$\tilde{\Phi}(s)=s\ln(\ln(s+\e)).$$ It is easy to check that $\tilde{\Phi}$ is a Young function. 
 Define the sequence $k=(k_n)_{n \in \mathbb{N}}$ by $k_n=\frac{\ln(Cn)}{n}$, $n\in\mathbb{N}$, where $C=\ln(2)+\ln(2\e)>1$. We choose $n$ large enough, such that $k_n n = \ln(Cn)\geq 1$ holds.  
 Similar to \cite[Ex.~4.2.13]{Diss.Wintermayr} one can show
 \begin{equation*}
 \tilde{\Phi}\left(\frac{2^n}{k_n n}\e^{-2^n t }\right) \leq 2^n\e^{-2^n t}.
 \end{equation*}
We deduce
\begin{align*}
\int_0^t \tilde{\Phi}\left( \frac{\e^{-2^n (t-s)}\frac{2^n}{n}}{k_n} \right) \, \mathrm{d}s \leq 1-\e^{-2^n t} < 1
\end{align*}
and hence $\|\e^{-2^n (t-\cdot)}\frac{2^n}{n}\|_{L_{\tilde{\Phi}}(0,t;\mathbb{C})} \leq k_n$ for sufficiently large $n$. Using the generalized H\"older inequality \eqref{Holder}, we get for $u \in E_\Phi(0,t;\ell^2(\mathbb{N}))$ and sufficiently large $n$
\begin{align*}
\left| \left(\int_0^t T_{-1}(t-s) B u(s) \, \mathrm{d}s \right) (n) \right|
&= \left|\int_0^t \e^{-2^n (t-s)} \frac{2^n}{n} (u(s)) (n) \, \mathrm{d}s \right|\\
&\leq 2 \left\|\e^{-2^n (t-\cdot)}\frac{2^n}{n}\right\|_{L_{\tilde{\Phi}}(0,t;\mathbb{C})} \|(u(\cdot))(n)\|_{E_\Phi(0,t;\mathbb{C})}  \\
&\leq 2 k_n \|u\|_{E_\Phi(0,t;\ell^2)},
\end{align*}
where we used in the last inequality that
\[ \int_0^t \Phi\left( \frac{|(u(s))(n)|  }{k} \right) \, \mathrm{d}s \leq \int_0^t \Phi\left( \frac{\|u(s)\|_{\ell^2}}{k} \right) \, \mathrm{d}s. \]
Therefore, for some $M>0$,
\[ \left\| \int_0^t T_{-1}(t-s) B u(s) \, \mathrm{d}s \right\|_{\ell^2} \leq M \|k\|_{\ell^2} \|u\|_{E_\Phi(0,t;\ell^2)}, \]
which shows that $B$ is $E_\Phi$-admissible and thus $\Sigma(A,[B,0],F)$ is integral ISS.
\end{example}

\section{Controlled Fokker--Planck equation}\label{sec:FP}
Following \cite{BreitenKunisch,Hartmannetal} we consider the following  variant of the Fokker--Planck equation on a bounded domain $\Omega\subset \mathbb R^n$, with smooth boundary $\partial\Omega$, 
\begin{equation}\label{eq:FokkerPlanck}
\arraycolsep=1.3pt\def\arraystretch{1.2}
\begin{array}{llll}
\frac{\partial \rho}{\partial t}(x,t) &=&{} \nu\Delta\rho(x,t)+\nabla\cdot\Big(\rho(x,t)\nabla V(x,t)\Big),\\
\rho(x,0) &=&{}  \rho_0(x), 
\end{array}
\end{equation}
where $ x\in \Omega, t>0$,   
with reflective boundary conditions 
\begin{equation}\label{eq:BCFP}  0=(\nu\nabla\rho  +\rho\nabla V)\cdot \vec{n}, \qquad \text{on }\partial\Omega\times(0,\infty). 
\end{equation}
Here $\vec{n}$ refers to the outward normal vector on the boundary, $\rho_0$ denotes the initial probability distribution with $\int_\Omega \rho_0(x) \, \mathrm{d}x=1$ and $\nu>0$. Furthermore, the potential $V$ is assumed to be of the form
 \begin{equation}\label{eq:potential}
 V(x,t)=W(x)+\alpha(x)u(t),
 \end{equation}
 where $W \in W^{2,\infty}(\Omega)$ and $\alpha \in W^{1,\infty} \cap H^2(\Omega)$ satisfying the structural assumption $\nabla \alpha \cdot \vec{n} =0$ on $\partial \Omega$. Thus, the scalar-valued input function $u$ enters via the spatial profile $\alpha$ in the potential.  
In order to cast the equations in an abstract framework, we introduce the following operators:
\begin{align*} 
{A}f={}& \nu\Delta f+\nabla\cdot(f\nabla W), \\
D({A})={}& \{f\in H^{1}(\Omega)\mid \Delta f\in L^{2}(\Omega), (\nu \nabla f+f \nabla W)\cdot \vec{n} =0 \text{ on } \partial \Omega\},\\
{Bf}={}& \nabla \cdot (f \nabla \alpha ), \\ 
D({B})={}& H^{1}(\Omega),
\end{align*}
where ${X}=L^{2}(\Omega)$ and $H^{1}(\Omega)$, $H^{2}(\Omega)$ refer to the standard Sobolev spaces. 
By standard arguments, the operator $A$ is seen to generate a bounded $C_0$-semigroup on $X$, with discrete spectrum $\sigma(A) = \sigma_p(A) \subseteq (-\infty,0]$ and $\rho_\infty = c\e^{-\Phi}$ is an eigenfunction to the simple eigenvalue $0$, where $\Phi$ is given by $\Phi = \ln \nu + \frac{W}{\nu}$ and $c>0$ is such that $\int_{\Omega}\rho_{\infty}\mathrm{d}x=1$, see \cite{BreitenKunisch}. Furthermore, we will identify $B$ with its extension from $X$ to $X_{-1}$. 

We now consider the system 
around the stationary distribution $\rho_{\infty}$ instead of the origin,  see also \cite{BreitenKunisch},
and decompose $X$ according to the projections
\[P:L^2(\Omega) \rightarrow L^2(\Omega), y\mapsto y -  \int_\Omega y(x) \, \mathrm{d}x \rho_\infty \quad \text{ and }\qquad Q:= I-P.\]
Note that $\ran(Q)=\ker(P)=\mathrm{span}{\{\rho_\infty\}}$ and $\ker(Q) = \ran(P)$. Define $ \mathcal{X}=\ran(P)$.
Using 
 $y=\rho - \rho_\infty$ and $y=y_P + y_Q$ with $y_P=Py \in \mathcal{X}$ and $y_Q=Qy\in\mathrm{span}{\{\rho_\infty\}}$ and following \cite[Sec.~3.2]{BreitenKunisch},  
the Fokker--Planck equation can be rewritten as
\begin{equation}\label{eq:FKfinal}
\arraycolsep=1.3pt\def\arraystretch{1.2}
\begin{array}{lll}
\dot{y}_P(t) &=\mathcal{A}y_P(t)+\mathcal{B}_1\left(y_P(t) u(t)\right) + \mathcal{B}_2u, \ t\geq0,\\
 y_P(0)&= P\rho_0,\\
 y_Q(t) &= Q\rho_0 - \rho_\infty = 0, \ t\geq0, 
 \end{array}\hspace{-0.3cm}
\end{equation}
where 
\begin{align*}
	&\mathcal{A}:D(\mathcal{A})=\mathcal{X}\cap D(A)\to \mathcal{X}, f\mapsto Af,\\
	&\mathcal{B}_{1}:\mathcal{X}\to \mathcal{X}_{-1}, f\mapsto Bf,
	\\ &\mathcal{B}_{2}:\C\to \mathcal{X}, u\mapsto uB\rho_{\infty},
\end{align*}
and $\mathcal{A}$ generates a strongly continuous semigroup as $APf=PAf$ for $f\in D(A)$, see also \cite[Eq.~(3.12)]{BreitenKunisch}. We emphasize that $Q\rho_0 - \rho_\infty = 0$ follows by the assumption that $\int_{\Omega}\rho_{0}(x)\mathrm{d}x=1$. That $\mathcal{B}_{1}$ and $\mathcal{B}_{2}$ are well-defined will be argued below.   

\begin{theorem} \label{thm:FP}
System \eqref{eq:FKfinal} is integral ISS.
\end{theorem}
In the remainder of the section we will lay out the proof of Theorem \ref{thm:FP} based on Theorem \ref{thm2} and Proposition \ref{prop:linear}. This includes to show that $\mathcal{A}$, $\mathcal{B}_{1}$, $\mathcal{B}_{2}$ satisfy the the assumptions of the abstract system class from Section \ref{sec2.1} and considering the nonlinearity 
\begin{equation}\label{eq:FFP}
	F:\mathcal{X}\times \mathbb{C} \to \mathcal{X}, (y,u)\mapsto  y u.
\end{equation} 

Let $M$ be the multiplication operator by $\e^{\frac{\Phi}{2}}$ on $L^2(\Omega)$. Clearly, $M$ is bounded and invertible on $L^2(\Omega)$, leaves $H^1(\Omega)$ invariant, and the inverse $M^{-1}$ is the multiplication operator by $\e^{-\frac{\Phi}{2}}$. Hence, $\tilde{A}$ given by
\begin{align*}
\tilde{A} =&{} M A M^{-1},\\
D(\tilde{A}) =&{} MD(A)
\end{align*}
is well-defined and self-adjoint.

To study admissibility of $B$ we introduce the following well-known abstract interpolation and extrapolation spaces, see e.g.\ \cite{TucWei09}.
Let $\tilde{X}_{1}$ and $\tilde{X}_{-1}$ be defined in the same way as $X_1$ and $X_{-1}$, but using $\tilde{A}$ instead of $A$. 
We define $\tilde{X}_{\frac{1}{2}}$ as the completion of $D(\tilde{A})$ with respect to the norm given by
\[\| z \|_{\tilde{X}_{\frac{1}{2}}}^2 := \langle (I-\tilde{A})z,z  \rangle,\quad x \in D(\tilde{A}),\] 
and we denote by $\tilde{X}_{-\frac{1}{2}}$ the dual space of $\tilde{X}_{\frac{1}{2}}$ with respect to the pivot space $X$, i.e. the completion of $X$ with respect to the norm
$ \sup_{\|v\|_{\tilde{X}_{\frac{1}{2}}}\leq 1} |\langle z,v \rangle_X|$. 
 The following embeddings are dense and continuous:
$\tilde{X}_1 \hookrightarrow \tilde{X}_{\frac{1}{2}} \hookrightarrow X \hookrightarrow \tilde{X}_{-\frac{1}{2}} \hookrightarrow \tilde{X}_{-1}$.

We first prove that the operator $\tilde{B}:= M B M^{-1}$ defined on $D(\tilde{A})$ has a unique extension $\tilde{B} \in L(X,\tilde{X}_{-\frac{1}{2}})$ which is $L^{2}$-admissible for $\tilde A$. Integration by parts gives
\[ \|v\|_{\tilde{X}_{\frac{1}{2}}}^2 = \|v\|_{L^2}^2 + \|\nabla\left(\e^{\frac{\Phi}{2}}v\right)e^{-\frac{\Phi}{2}}\|_{L^2}^2, \qquad v\in D(\tilde{A}). \]
For $f \in D(\tilde{A})$ and $v \in D(\tilde{A})$, $\| v\|_{\tilde{X}_{\frac{1}{2}}} \leq 1$, we have that
\begin{eqnarray*}
| \langle \tilde{B} f,v \rangle_{L^2} | &=&
\left|\int_\Omega v \e^{\frac{\Phi}{2}} \nabla \cdot \left( \e^{-\frac{\Phi}{2}} f \nabla \alpha \right) \, \mathrm{d}x \right| \\
&=& \left|\int_{\partial \Omega} v \e^{\frac{\Phi}{2}} \e^{-\frac{\Phi}{2}} f \nabla \alpha \cdot \vec{n} \, \mathrm{d}\sigma - \int_\Omega \nabla\left( v \e^{\frac{\Phi}{2}}\right)\cdot \left(\e^{-\frac{\Phi}{2}} f \right) \, \mathrm{d}x \right|\\
&\leq& \| \nabla\left( v \e^{\frac{\Phi}{2}} \right) \e^{-\frac{\Phi}{2}} \|_{L^2(\Omega)^n}^2 \|f \nabla \alpha\|_{L^2(\Omega)^{n}}^2 \\
& \leq& n \| \nabla \alpha\|_{L^2(\Omega)^n}^2 \| \nabla\left( v \e^{\frac{\Phi}{2}} \right) \e^{-\frac{\Phi}{2}} \|_{L^2(\Omega)^n}^2 \|f \|_{L^2(\Omega)}^2
\end{eqnarray*} where $\sigma$ is the surface measure on $\partial \Omega$. Thus $\tilde{B}\in L(X,\tilde{X}_{-\frac{1}{2}})$ and  $\tilde{B}$ is $L^2$-admissible for $\tilde{A}$  by \cite[Prop.~5.1.3]{TucWei09}. We have for $\beta \in \rho(A)=\rho(\tilde{A})$ and $f \in X$
\begin{align*}
\|M^{-1}f\|_{X_{-1}} &= \| (\beta-A)^{-1} M^{-1}f\|_X = \| M^{-1}(\beta-\tilde{A})^{-1}f\|_X \leq\|M^{-1}\|\|f\|_{\tilde{X}_{-1}}.
\end{align*}
 Thus, $M^{-1}$ extends uniquely to an operator in $L(\tilde{X}_{-1},X_{-1})$. The same argument yields a unique extension  $M\in L(X_{-1},\tilde{X}_{-1})$. Note that these extensions are inverse to each other, so it is natural to denote the extensions again by $M$ and $M^{-1}$.\\
We claim that $M^{-1} \tilde{B} M \in L(X,X_{-1})$ extends $B$ to an $L^{2}$-admissible operator for $A$ which we again denote by $B$. 
 Indeed, if $(T(t))_{t \geq 0}$ is the semigroup generated by $A$, then $(S(t))_{t \geq 0}$ with $S(t)=MT(t)M^{-1}$ is the semigroup generated by $\tilde{A}$ and for $u \in L^2(0,t;X)$ we have $ M u \in L^2(0,t;X)$ and
\[\int_0^t T_{-1}(t-s)B u (s) \, \mathrm{d}s = M^{-1} \int_0^t S(t-s) \tilde{B}(Mu)(s) \, \mathrm{d}s.\]


As $\mathcal{B}_2 \in L(\mathbb{C},\mathcal{X})$, $\mathcal{B}_2$ is clearly $L^1$-admissible.
The operator $P$ commutes with the $C_0$-semigroup generates by $A$ \cite[Eq.~(3.12)]{BreitenKunisch}, by  \cite[Lem.~4.4]{JacobZwart} the operator $\mathcal{B}_1=B|_{\mathcal{X}}\in \mathcal{L}(\mathcal{X},\mathcal{X}_{-1})$ is well-defined and $L^2$-admissible for $\mathcal{A}$.

Thus the bilinearly controlled Fokker--Planck system given by \eqref{eq:FokkerPlanck}-\eqref{eq:potential} can be written as a system $\Sigma(\mathcal{A},[\mathcal{B}_{1},\mathcal{B}_{2}],F)$. 

Remark \ref{rem:aftermainthm}  implies that the Fokker--Planck system \eqref{eq:FokkerPlanck}-\eqref{eq:potential} has a unique global mild solution $\rho$ for any initial value $\rho_{0}\in L^{2}(\Omega)$ and input function $u\in L^{2}(0,\infty;U)$. Further, in \cite[Proposition 2.2]{BreitenKunisch} it is shown that $\int_{\Omega}\rho_{0}(x)\mathrm{d}x=1$ implies $\int_{\Omega}\rho(t,x)\mathrm{d}x=1$ for all $t>0$.

Following the construction of the integral ISS estimate (c.f.~\eqref{bilinearISSestimate}) we deduce an explicit integral ISS estimate:
There exists constants $C,\omega>0$ such that for any $\rho_{0}\in L^{2}(\Omega)$ with $\int_{\Omega}\rho_{0}(x)\mathrm{d}x=1$ and $u\in L^2(0,\infty;U)$, the 
 global mild solution of the Fokker--Planck system \eqref{eq:FokkerPlanck}  satisfies
\begin{align*}
\|\rho(t)-\rho_{\infty}\|_{L^{2}}\le& \ C \e^{-\omega t} \left(\|\rho_0-\rho_{\infty}\|_{L^{2}} + \|\rho_0-\rho_{\infty}\|_{L^{2}}^{2}\right)+\gamma\left(\int_{0}^{t}\|u(s)\|_{U}^{2}\mathrm{d}s\right),
 \end{align*}
 where
$ \gamma(r) =C r\e^{Cr^{\frac{1}{2}}}+Cr^{\frac{1}{2}}+Cr$.

\section{Conclusion}
Bilinear systems appear naturally in control theory e.g.\ when considering multiplicative disturbances in feedback loops of linear systems.
The results in this article draw a link between bilinear systems, which are a classical example class in (integral) ISS in finite-dimensions, and recent progress in ISS for infinite-dimensional systems. We emphasize that the most natural example in this context, 
\[\dot{x}(t)=Ax(t)+u(t)x(t),\quad t>0,\quad x(0)=x_{0},\]
with $A$ generating a $C_{0}$-semigroup $(T(t))_{t \geq 0}$ on $X$, is covered by the system class considered here. More precisely, by the results in Section \ref{sec2}, it follows that this system is integral ISS if and only if $(T(t))_{t \geq 0}$ is exponentially stable. More precisely, the sufficiency follows since the identity is $L^1$-admissible and hence the system is integral ISS. It seems that prior works on integral ISS \cite[Sec.~4.2]{MiI14b} did not cover this comparably simple class as the bilinearity  $x\mapsto xu$ fails to satisfy a Lipschitz condition uniform in $u$ required there\footnote{However, it seems that this can be overcome with a  carefully refined argument in the proof of \cite[Thm.~4.2]{MiI14b}.}.

Moreover,  our results generalize to integral ISS assessment for bilinearities arising from boundary control (or lumped control).

\section{Appendix}\label{appendix}

We briefly introduce Orlicz spaces of functions $f:I\to Y$ for an interval $I\subset \mathbb{R}$ and a Banach space $Y$. For more details on Orlicz spaces we refer to \cite{Adams,KrasnRut,Kufner}.
Let $\Phi:\Rpn \rightarrow \Rpn$ be a {\em Young function}, i.e.\ $\Phi$ is continuous, increasing, convex with $\lim_{s \to 0}\frac{\Phi(s)}{s}=0$ and $\lim_{s \to \infty}\frac{\Phi(s)}{s}= \infty$ and denote by $L_\Phi(I;Y)$ the set of Bochner-measurable functions $u:I\rightarrow Y$ for which there exists a constant $k>0$ such that  $\Phi(k\|u(\cdot)\|)$ is integrable. We equip $L_\Phi(I;Y)$ with the norm
\begin{equation}\label{defOrlicz}
 \|u\|_{L_\Phi(I;Y)} = \inf \left\{ k>0 \middle| \int_I \Phi\left(\frac{\|u(s)\|}{k}\right)\, \mathrm{d}s \leq 1 \right\}.
 \end{equation}
Despite the fact that $L_\Phi(I;Y)$ is typically referred to as ``Orlicz space'' in the literature, we prefer to call
\[E_\Phi(I;Y) = \overline{\{u \in L^\infty(I;Y) \mid \esssupp\, u \text{ is bounded} \}}^{\|\cdot\|_{L_\Phi(I;Y)}}\]
the {\em Orlicz space} associated with the Young function $\Phi$. We write $\|u\|_{E_\Phi(I;Y)} = \|u\|_{L_\Phi(I;Y)}$  for $u \in E_\Phi$. Note that $u \in E_\Phi(I;Y)$ implies that $\Phi \circ \|u(\cdot)\|$ is integrable.  Typical examples of Orlicz spaces are $L^{p}$-spaces; for $\Phi(t)=t^{p}$ with $p\in(1,\infty)$ it holds that $E_{\Phi}(I;Y)$ is isomorphic to $L^{p}(I;Y)$. 
\\
A Young function $\Phi$ is said to satisfy the {\em $\Delta_2$-condition} if there exist $K>0$ and $s_{0}\ge 0$ such that \[\Phi(2s) \leq K \Phi(s),\qquad s\geq s_0.\] 
Note that $E_\Phi(I;Y) = L_\Phi(I;Y)$ if and only if $\Phi$ satisfies the $\Delta_2$-condition. Also note that 
$\Phi(s)=s^p$, $p \in (1,\infty)$ satisfies the $\Delta_{2}$-condition.
For a Young function $\Phi$ the {\em complementary Young function} $\tilde{\Phi}$ is defined by
$$\tilde{\Phi}(s)= \max_{t \geq 0}\left( st - \Phi(t) \right).$$
Again, this is a Young function and $\Phi$ can be recovered from $\tilde{\Phi}$ in the same manner. The complementary Young function to $\Phi(s)=\frac{s^p}{p}$, $1<p<\infty$, is given by $\tilde{\Phi}(s)=\frac{s^q}{q}$ with $\frac{1}{p}+ \frac{1}{q}=1$.\\
As for $L^p$ spaces, an equivalent norm to $\|\cdot\|_{L_\Phi(I;Y)}$ is given by
\begin{equation}\label{aequinorm}
 \| u \|_{\Phi,(I;Y)} = \sup\left\{ \int_I \|u(s)\| |v(s)| \, \mathrm{d}s \middle| v \text{ measurable, } \int_I \tilde{\Phi}(|v(s)|) \, \mathrm{d}s \leq 1 \right\}. 
\end{equation}
Furthermore, for a  Young functions $\Phi$ and its complementary Young function $\tilde{\Phi}$ the following generalized H\"older inequality
\begin{equation}\label{Holder}
\int_I \|u(s)\| \|v(s)\| \mathrm{d}s \leq 2 \|u\|_{L_\Phi} \|v\|_{L_{\tilde{\Phi}}}. 
\end{equation}
holds. This also implies the continuity of the embeddings 
\[L^{\infty}(I;Y)\hookrightarrow L_\Phi(I;Y) \hookrightarrow L^1(I;Y)\] if $I$ is bounded. Although $L^{1}$ is not an Orlicz space, we will explicitly allow for $\Phi(t)=t$ in our notation referring to $E_\Phi(I;Y)=L^1(I;Y)$. Note that the definition of the norm \eqref{defOrlicz} is indeed consistent with the $L^1$-norm and that $\Phi$ satisfies the $\Delta_2$-condition. However, we will not define a ``complementary Young function" for this particular $\Phi$. \\
An essential property of Orlicz spaces is the absolute continuity of the $E_\Phi$ norm with respect to the length of the interval $I$ (see e.g. \cite[Thm.~3.15.6]{Kufner}), this is for $u \in E_\Phi(I;Y)$ and $\varepsilon > 0$ there exists $\delta > 0$ such that for each interval $I$ holds
\[ \lambda(I) < \delta \quad \implies \quad \|u\|_{E_\Phi(I;Y)} < \varepsilon,\]
where $\lambda$ refers to the Lebesgue-measure on $\mathbb{R}$.

\def\cprime{$'$}



\end{document}